\documentclass[a4paper,12pt]{article}

\usepackage[top=3.0cm,bottom=3.0cm,left=2.5cm,right=2.5cm]{geometry}
\usepackage{amsmath,amsthm,mathrsfs,graphicx,amsfonts}
\usepackage{bm}
\usepackage{cases}
\usepackage[english]{babel}
\usepackage{amsmath,amsthm}
\usepackage{amsfonts}
\usepackage{latexsym}
\usepackage{graphicx}
\usepackage{txfonts}
\usepackage[numbers,sort&compress]{natbib}
\usepackage[natural]{xcolor}
\usepackage{rotating}
\usepackage{mathtools}

\newtheorem{theorem}{Theorem}[section]
\newtheorem{corollary}[theorem]{Corollary}
\newtheorem{lemma}[theorem]{Lemma}
\newtheorem{claim}{}[theorem]
\newtheorem{conjecture}[theorem]{Conjecture}

\newcommand{\del}{\backslash}

\newcommand{\cB}{\mathcal{B}}

\newcommand{\cl}{\hbox{\rm cl}}

\newcommand{\st}{\hbox{\rm st}}

\newcommand{\cC}{\mathcal{C}}
\newcommand{\cI}{\mathcal{I}}

\newcommand{\loops}{\hbox{\rm loops}}

\begin{document}

\title{The $9$-connected Excluded Minors for the Class of Quasi-graphic Matroids}

\author{Rong Chen
\footnote{Center for Discrete Mathematics, Fuzhou University, Fujian 350003, China.
Email: rongchen@fzu.edu.cn.
This research was partially supported by grants from the National Natural Sciences
Foundation of China (No. 11971111) and NSFFP (No. 2019J01645).} }

\date{}
\maketitle
\begin{abstract}
The class of quasi-graphic matroids, recently introduced by Geelen, Gerards, and Whittle, is minor closed
and contains both the class of lifted-graphic matroids and the class of frame matroids, each of which
generalises the class of graphic matroids. In this paper, we prove that the matroids $U_{3,7}$ and $U_{4,7}$ are
the only $9$-connected excluded minors for the class of quasi-graphic matroids.
\end{abstract}



\section{Introduction}

Let $H$ be a graph and let $N$ be a matroid. For a vertex $v$ of $H$ we let $\loops_H(v)$
denote the set of loops of $H$ whose ends are $v$. We say that $H$ is a
{\em framework} for $N$ if
\begin{itemize}
\item[(QG1)] $E(H)=E(N)$,
\item[(QG2)] $r_N(E(H'))\leq |V(H')|$ for each component $H'$ of $H$, and
\item[(QG3)] for each vertex $v$ of $H$ we have
$\cl_N(E(H-v))\subseteq  E(H-v)\cup {\rm loops}_H(v)$, and
\item[(QG4)] for each circuit $C$ of $N$, the graph $H[C]$ has at most two components.
\end{itemize}
A matroid is {\em quasi-graphic} if it has a framework. The class of quasi-graphic matroids, recently introduced
by Geelen, Gerards, and Whittle ~\cite{GGW}, is minor closed and contains both lifted-graphic matroids and frame
matroids. Recently, the author and Geleen \cite{CG17} proved that there are infinitely many quasi-graphic excluded
minors for the class of frame matroids and the class of lifted-graphic matroids, but we are confident that the
class of quasi-graphic matroids admits a finite excluded-minor characterisation.

\begin{conjecture}(\cite{CG17}, Conjecture 1.5.)\label{conj-quasi}
There are, up to isomorphism, only finitely many excluded-minors for the class of quasi-graphic matroids.
\end{conjecture}

One of the difficulties to prove Conjecture \ref{conj-quasi} is that some graphic matroids have exponentially many
different frameworks; for example, the rank-$r$ wheel has at least $2^r$ ``inequivalent" frameworks, see \cite{CDFP}. The
same difficulty appears when considering problems on excluded minors for the class of frame matroids and for the
class of lifted-graphic matroids. In fact, in the proof of Rota's Conjecture, Geleen, Gerards, and Whittle encountered
a similar difficulty. The interesting thing is: we have some kind of opposite versions in the proof of the two conjectures.
For Rota's Conjecture, the proof for the low branch-width case is not complicated, see \cite{BG16,GW02}; while the proof
for the high branch-width case is very difficult. While, for Conjecture \ref{conj-quasi}, the proof for low connectivity
is thought to be difficult, while the proof for high connectivity is not complicated. In this paper, we prove

\begin{theorem}\label{main-result}
Other than $U_{3,7}$ and $U_{4,7}$, no excluded minor for the class of quasi-graphic matroids is $9$-connected.
\end{theorem}

Funk and Mayhew  \cite{FM17} recently proved that, for each positive integer $r$, the class of quasi-graphic
matroids has only a finite number of excluded minors of rank $r$.

This paper is organized as follows. In Section 3, we prove that $U_{3,7}$ and $U_{4,7}$ are the only 9-connected
excluded minors of rank less than nine for the class of quasi-graphic matroids. 9-connected excluded minors of
rank at least nine are considered in Section 5. Some definitions and basic properties of quasi-graphic matroids
are given in Section 2. Properties of frameworks for graphic matroids are presented in Section 4.

\section{Preliminaries}
We assume that the reader is familiar with matroid theory and we follow the terminology of Oxley~\cite{Oxley}.

For a graph $G$, let $\loops(G)$ be the set of loops in $G$. An edge of $G$ is a {\sl link} if it is not a loop. For
any $v\in V(G)$,  let $\st_G(v)$ denote the set of edges incident with $v$.
For any $U\subseteq V(G)$ and $F\subseteq E(G)$, set $\st_G(U)=\bigcup_{u\in U}\st_G(u)$, and let $G[U]$ be the induced
subgraph of $G$ defined on $U$, and let $G[F]$ be the subgraph of $G$ with $F$ as its edge set and without isolated vertices.
Let $c_G(F)$ be the number of components of $G[F]$, and let $V_G(F)$ denote $V(G[F])$.
When $F=\{e\}$, we will let $V_G(e)$ denote $V_G(\{e\})$. When there is no confusion, all subscripts will be omitted.
For a number $k$, we say that $G$ is $k$-connected if $G-S$ has exactly one component for any $S\subset V(G)$ with $|S|<k$.

A {\em theta graph} is a graph that consists of a pair of distinct vertices joined by three internally disjoint
paths. A {\em cycle} is a connected $2$-regular graph. A collection ${\mathcal B}$ of some cycles of $G$ satisfies
the {\em theta property} if no theta subgraph of $G$ contains exactly two members of $\mathcal B$.
A \emph{biased graph} consists of a pair $(G, \mathcal{B})$, where
$G$ is a graph and $\mathcal{B}$ is a collection of some cycles of $G$ that satisfies the theta
property. A cycle $C$ of $G$ is {\em balanced} if $C\in\cB$, otherwise, it is {\em unbalanced}.


Let $H$ be a framework for a matroid $N$. For any cycle $C$ of $H$, either $C\in\cC(N)$ or $C\in\cI(N)$ by (\cite{GGW}, Lemma 2.5.).
Let $\cB_N$ be the set of cycles of $H$ that are circuits of $N$. Since $\cB_N$ satisfies the theta property
by (\cite{GGW}, Lemma 3.2.), $(H,\cB_N)$ is a biased graph. For convenience, we will also view $H$ as the biased
graph $(H,\cB_N)$. A subgraph $H'$ of $H$ is {\em balanced} if each cycle in $H'$ is balanced; otherwise, $H'$ is {\em unbalanced}.
If all cycles in $H'$ are unbalanced, then $H'$ is {\em contra-balanced}.

By (\cite{GGW}, Lemma 3.3) and (QG4), we have

\begin{lemma}\label{circuit in framework}
Let $H$ be a framework for a matroid $N$. When $C\in\cC(N)$, either

\begin{enumerate}
\item $H[C]$ is a balanced cycle,
\item $H[C]$ is a connected contra-balanced graph with minimum degree at least two with $|C|=|V(C)|+1$, or
\item $H[C]$ is a union of two unbalanced cycles that meet in at most one single vertex.
\end{enumerate}
\end{lemma}

\begin{lemma}(\cite{GGW}, Lemma 2.6.)\label{find circuit}
Let $H$ be a framework for a matroid $N$. If $H'$ is a subgraph of $H$ with $|E(H')|>|V(H')|$, then $E(H')$ is a dependent set of $N$.
\end{lemma}

By Lemmas \ref{circuit in framework} and \ref{find circuit}, we have

\begin{lemma}\label{determine circuit}
Let $H$ be a framework for a matroid $N$. Let $C_1, C_2$ be unbalanced cycles of $H$ with $|V(C_1)\cap V(C_2)|\leq1$. Then the following hold.
\begin{itemize}
\item $E(C_1\cup C_2)$ is a circuit of $N$ when $|V(C_1)\cap V(C_2)|=1$.
\item When $|V(C_1)\cap V(C_2)|=0$, for each minimal path $P$ in $H$ linking $C_1$ and $C_2$, we have $E(C_1\cup C_2)\in\cC(N)$ or $E(C_1\cup C_2\cup P)\in\cC(N)$.
\end{itemize}
\end{lemma}

We say that $H$ is a {\em frame representation} of a matroid $N$ if a subset $I$ of $E(H)$ is independent in $N$ if and
only if $H[I]$ has no balanced cycles and $|E(H')|\leq|V(H)|$ for each component $H'$ of $H[I]$. We say that $H$ is a {\em lifted-graphic
representation} of $N$ if a subset $I$ of $E(H)$ is independent in $N$ if and only if $H[I]$ has at most one cycle and when
the cycle exists, it is unbalanced. Note that, when $H$ is a lifted-graphic representation for a 3-connected matroid, $H$ has at most one loop, and the loop is unbalanced.

\begin{theorem}(\cite{GGW}, Theorems 7.1 and 7.2.)\label{ub loop}
Let $H$ be a framework for a $3$-connected matroid $N$. If $H$ has an unbalanced loop, then $H$ is a frame representation or a lifted-graphic representation for $N$.
\end{theorem}

When $H$ is a lifted-graphic representation for $N$ with an unbalanced loop $e$, by the definition of lifted-graphic representation, all graphs obtained from $H\del e$ by attaching the loop $e$ to any vertex of $H\del e$ or a new vertex not in $H\del e$ are also lifted-graphic representations of $N$. Under this condition, we view all graphs obtained in this way as equivalent. That is, when all frameworks for $N$ can be obtained from $H$ by this way, we view $H$ as the  unique framework for $N$.

\begin{lemma}\label{disconnected}
Let $H$ be a framework for a matroid $N$. If $H$ is not connected but $N$ is connected, then $H$ is a lifted-graphic representation of $N$.
\end{lemma}

\begin{proof}
Let $H_1,H_2,\ldots, H_n$ be the components of $H$. Since every pair of elements of $E(N)$ must be contained in a circuit of $N$,  by Lemma \ref{circuit in framework} each edge of $H$ is in an unbalanced cycle. 

\begin{claim}\label{2 component}
Let $C_1$ and $C_2$ be unbalanced cycles of $H_1$ and $H_2$, respectively. If $H_2$ has an unbalanced cycle $C'_2$ with $E(C_2)\cap E(C'_2)\neq\emptyset$ and $E(C_1\cup C'_2)\in\cC(N)$, then $E(C_1\cup C_2)\in\cC(N)$.
\end{claim}
\begin{proof}[Subproof.]
Assume not. Without loss of generality we may assume that $C'_2$ is chosen with $E(C_2\cup C'_2)$ as small as possible. When  $C_2\cup C'_2$ is a theta subgraph of $H_2$,  since $E(C_2\cup C'_2)$ or the third cycle in $C_2\cup C'_2$ that is neither $C_2$ nor $C'_2$ is a circuit of $N$ by Lemmas \ref{circuit in framework} and \ref{find circuit}, we have $E(C_1\cup C_2)\in\cC(N)$ by the circuit elimination axiom and Lemma \ref{circuit in framework}. So we may assume that $C_2\cup C'_2$ is not a theta subgraph of $H_2$. Since $E(C_2)\cap E(C'_2)\neq\emptyset$, there is a path $P\subsetneq C_2$ such that $C'_2\cup P$ is a theta-graph. In a similar way we can show that $E(C_1\cup C_2^{''})$ is a circuit of $N$ for an unbalanced cycle $C_2^{''}$ of $H_2$ with $P\subseteq C_2^{''}\subseteq C'_2\cup P$, a contradiction to the choice of $C'_2$ as $|E(C_2\cup C_2^{''})|<|E(C_2\cup C'_2)|$.
\end{proof}

\begin{claim}\label{circuit}
A union of each pair of unbalanced cycles coming from different components of $H$ is a circuit of $N$.
\end{claim}
\begin{proof}[Subproof.]
Let $C_1$ and $C_2$ be unbalanced cycles of $H_1$ and $H_2$, respectively. By symmetry,  it suffices to show that  $E(C_1\cup C_2)\in\cC(N)$. Let $e_i\in C_i$ for each $1\leq i\leq2$. Since $N$ is connected, $N$ has a circuit $C$ containing $\{e_1,e_2\}$. Since $e_1$ and $e_2$ are in different components of $H$, by Lemma \ref{circuit in framework}, there is an unbalanced cycle $C'_i$ of $H_i$ containing $e_i$ for each integer $1\leq i\leq2$ such that $C=E(C'_1\cup C'_2)$. Since $e_2\in E(C_2)\cap E(C'_2)$, we have $E(C'_1\cup C_2)\in\cC(N)$ by \ref{2 component}. Moreover, since $e_1\in E(C_1)\cap E(C'_1)$,  using \ref{2 component} again, $E(C_1\cup C_2)\in\cC(N)$.
\end{proof}

\begin{claim}\label{circuit+1}
For every $1\leq i\leq n$, a union of every pair of vertex-disjoint unbalanced cycles of $H_i$ is a circuit of $N$.
\end{claim}
\begin{proof}[Subproof.]
Assume that the claim does not hold for $H_1$. Then there are vertex-disjoint unbalanced cycles $C_1, C'_1$ of $H_1$ and a path $P$ minimal linking the two cycles such that $E(C_1\cup C'_1\cup P)$ is a circuit of $N$ by Lemma \ref{determine circuit}. Let $C$ be a union of $C_1$ and an unbalanced cycle of $H_2$. By \ref{circuit}, $C$ is a circuit of $N$. Let  $f\in E(C_1)$ and $g\in E(P)$. By circuit elimination axiom, there is a circuit $C'$ of $N$ with $g\in C'\subseteq E(C_1\cup C'_1\cup P\cup C)-\{f\}$, a contradiction to Lemma \ref{circuit in framework} as $H[C']$ has degree-1 vertices.
\end{proof}

By \ref{circuit} and \ref{circuit+1}, a union of every pair of vertex-disjoint unbalanced cycles of $H$ is a circuit of $N$, so the lemma holds.
\end{proof}

After this paper was submitted to a journal in September 2017, one of the referees told the author in his/her referee report that Lemma \ref{disconnected} was also proved in (\cite{BFS}, Corollary 4.7) by  Bowler, Funk, and Slilaty. The two proofs are totally different.

By (\cite{GGW}, Lemmas 3.6 and 4.2) or Lemma \ref{disconnected} and (\cite{GGW}, Lemmas 4.2) we have

\begin{lemma}\label{3-con case}
Assume that $H$ is a framework for a $3$-connected matroid $N$ with $|E(N)|\geq4$ and $H$ has no isolated vertices. Then
\begin{enumerate}
\item $H$ is connected, or
\item $H$ is a lifted-graphic representation of $N$ with exactly two components, one of which is a loop-component.
\end{enumerate}
Moreover, $N$ has a connected framework.
\end{lemma}

\begin{lemma}\label{2-conn}
For any integer $k\geq 2$, if $H$ is a connected framework for a $k$-connected matroid $N$, then $H$ is $k-1$ connected.
\end{lemma}
\begin{proof}
Assume not. Let $(X,Y)$ be a partition of $E(N)$ with $m=|V_H(X)\cap V_H(Y)|\leq k-2$ and such that $H[X]$ and $H[Y]$ are connected graphs with at least $m+1$ vertices. When $H[X]$ and $H[Y]$ are unbalanced, implying that $|X|, |Y|\geq m+1$,  we have that $(X,Y)$ is an $m+1$-separation, a contradiction. When $X$ is balanced, $(X,Y)$ is an $m$-separation, a contradiction.
\end{proof}

\begin{theorem}(\cite{GGW}, Theorem 1.6.)\label{3-con case-1}
A $3$-connected matroid $N$ is quasi-graphic if and only if  there exists a graph $H$ such that
\begin{enumerate}
\item $E(H)=E(N)$,
\item $H$ is connected,
\item $r(N)\leq |V(H)|$, and
\item for each vertex $v$ of $H$ we have
$\cl_M(E(H-v))\subseteq  E(H-v)\cup {\rm loops}_H(v)$.
\end{enumerate}
\end{theorem}


\begin{lemma}\label{e in C}
Let $H$ be a framework for a matroid $N$. For an edge $e$ of $H$, if $H\del e$ is connected and unbalanced, then $e$ is in a circuit of $N$.
\end{lemma}
\begin{proof}
By considering a maximal independent set of $H\del e$, it follows from (QG2) that $r(N\del e)=|V(H)|$. Moreover, since $r(N)\leq |V(H)|$ by (QG2), we have $r(N)=r(N\del e)$. So the lemma holds.
\end{proof}

Let $H$ be a framework for a matroid $N$. Let $H'=H-\loops(H)$ when $H$ is a lifted-graphic representation of $N$, otherwise let $H'=H$. A vertex $v\in V(H)$ is a {\em blocking vertex} if $H'$ is unbalanced and all unbalanced cycles of $H'$ contain $v$. Set $\st_H^*(v)=\st_{H'}(v)$. Note that $\st_H^*(v)$ is the  same as $\st_H(v)$ unless $H$ is lifted-graphic of $N$ and $v$ is incident with a loop.


\begin{lemma}\label{star of blocking vertex}
Let $H$ be a connected framework for a $3$-connected matroid $N$, and $v\in V(H)$.
\begin{enumerate}
\item $\st^*_H(v)$ is a union of cocircuits of $N$.
\item$v$ is a blocking vertex of $H$ if and only if $\st^*_H(v)\notin\cC^*(N)$.
\end{enumerate}
\end{lemma}
\begin{proof}
(1) follows from Lemma \ref{circuit in framework} and Theorem \ref{ub loop}. Next, we  prove that (2) is true.

Note that $H$ and $H'$ are 2-connected by Lemma \ref{2-conn}.
Assume that $v$ is a blocking vertex of $H$. Since $H'-v$ is connected and balanced, $r(E(H'-v))=|V(H')|-2=r(N)-2$. So $\st^*_H(v)\notin\cC^*(N)$.

Assume that $\st^*_H(v)\notin\cC^*(N)$. Then $\st^*_H(v)$ contains at least two cocircuits of $N$ by (1), 
implying $r(E(H'-v))\leq r(N)-2=|V(H')|-2$. Moreover, since $H'-v$ is connected, $H'-v$ is balanced. So $v$ is a blocking vertex of $H$.
\end{proof}

Let $H$ be a connected framework for a 3-connected matroid $N$. We say that a vertex $v$ of $H$ is {\em fixed} in $H$ if $N\del \st^*_H(v)$ is a $3$-connected non-graphic matroid. 


\begin{lemma}\label{fix-0}
Let $H$ be a connected framework for a $3$-connected matroid $N$. For an edge $f$ of $H$, if $v$ is fixed in $H\del f$, then $v$ is fixed in $H$.
\end{lemma}
\begin{proof}
Evidently, it suffices to show that $N\del\st^*_H(v)$ is 3-connected. Assume not. Then $f\notin\st^*_H(v)$. Since $N\del (\st^*_H(v)\cup\{f\})$ and $N$ are 3-connected and non-graphic, $f$ is a coloop 
of $N\del\st^*_H(v)$ and $H'-\{v,f\}$ is connected and unbalanced by Lemma \ref{2-conn} and Theorem \ref{ub loop}. 
Then $f\in\st_H(v)$ by Lemma \ref{e in C}. Since $f\notin\st^*_H(v)$, we have that $\{f\}=\loops_H(v)$ and $H$ is a lifted-graphic representation for $N$. Hence, $f\in\cl(E(H-v))$  as $H-v$ is unbalanced, a contradiction to the fact that $f$ is a coloop of $N\del\st^*_H(v)$ .
\end{proof}

\begin{lemma}\label{fix}
Let $H$ and $H'$ be $2$-connected frameworks for a $3$-connected matroid $N$. If $v$ is a fixed vertex of $H$, then $\st_H^*(v)\in\cC^*(N)$ and there is a fixed vertex $v'$ of $H'$ satisfying $\st_{H'}^*(v')=\st_H^*(v)$.
\end{lemma}
\begin{proof}
Since $H'-v$ is unbalanced, $\st_H^*(v)\in\cC^*(N)$ by Lemma \ref{star of blocking vertex}.
So $r(N\del \st^*_H(v))=r(N)-1$. Since $N\del \st^*_H(v)$ is a $3$-connected non-graphic matroid, $|V(H')|=r(N)\geq3$ and by Lemma \ref{3-con case} the graph $H'\del \st^*_H(v)$ has exactly two components, one of which is an isolated vertex or a loop-component. Let $\{v'\}$ be the vertex set of the 1-vertex component of $H'\del \st^*_H(v)$, and $H'_1$ be the other component. Since no circuit of $N$ can intersect $\st_H^*(v)$ with exactly one element and $H'_1$ is connected and unbalanced, by Lemma \ref{e in C}, each edge in $\st_H^*(v)$ has at most one end in $H'_1$. So $\st^*_{H'}(v')=\st^*_H(v)$ by Theorem \ref{ub loop}, implying that $v'$ is fixed in $H'$.
\end{proof}

For convenience, we will say that the vertex $v'$ in Lemma \ref{fix} is the {\sl corresponding vertex} of $v$ in $H'$ and denote it by $v$ too.


\begin{lemma}\label{almost fix}
Let $H$ be a $2$-connected framework for a $3$-connected matroid $N$. If at most one vertex is not fixed in $H$, then one of the following holds.
\begin{enumerate}
\item $H$ is the unique framework for $N$.
\item $H-\loops(H)$ has a blocking vertex.
\end{enumerate}
In particular, when all vertices are fixed in $H$, (1) holds.
\end{lemma}
\begin{proof}
First, we prove
\begin{claim}\label{all fix}
When all vertices in $H$ are fixed, (1) holds.
\end{claim}
\begin{proof}[Subproof.]
When $\loops(H)=\emptyset$ or $H$ has a loop but it is not a lifted-graphic representation for $N$, Lemma \ref{fix} implies that (1) holds. So we may assume that $H$ is a lifted-graphic representation for $N$ with a loop $f$ by Theorem \ref{ub loop}. Let $H'$ be another connected framework for $N$. By Lemma \ref{fix}, $\loops(H)=\loops(H')=\{f\}$ and $H\del f=H'\del f$. When $H$ has no blocking vertex, since a union of $f$ and each unbalanced cycle of $H$ is a circuit of $N$ by the structure of $H$, the graph $H'$ must be a lifted-graphic representation for $N$, so $H$ and $H'$ are equivalent. When $H$ has a blocking vertex, either $f$ is incident with a blocking vertex of $H'\del f$ or $H'$ must be a lifted-graphic representation for $N$. No matter which case happens, $H'$ is a lifted-graphic representation for $N$, so $H$ and $H'$ are equivalent. That is, (1) holds.
\end{proof}

By \ref{all fix}, we may therefore assume that $H$ has a unique unfixed vertex $v$. Assume that (1) is not true. Let $H'$ be a connected framework for $N$ that is not equivalent to $H$.
By Lemma \ref{fix}, we may assume that $V(H)=V(H')$ and $\st^*_H(u)=\st^*_{H'}(u)$ for any $v\neq u\in V(H)$.
Therefore,

\begin{claim}\label{H-v}
For any vertices $x,y\in V(H-v)$, we have that $xy\in E(H)$ if and only if $xy\in E(H')$.
\end{claim}

By symmetry, \ref{all fix} and Lemma \ref{fix}, we may assume that $v$ is the unique unfixed vertex of $H'$.

\begin{claim}\label{1 case}
When $H$ is a lifted-graphic representation of $N$ with a loop $f$, both $H$ and $H'-\loops(H')$ have $v$ as their blocking vertex.
\end{claim}
\begin{proof}[Subproof.]
Since $N$ is 3-connected, $\loops(H)=\{f\}$. Since $\st_H(u)-\{f\}=\st^*_H(u)=\st^*_{H'}(u)$ for any $v\neq u\in V(H)$, we have $f\in\loops(H')$. When $H'$ is a lifted-graphic representation for $N$, since $\{f\}=\loops(H')$, we have $H\del f=H'\del f$, so $H$ and $H'$ are equivalent, a contradiction. Hence, $H'$ is a frame representation for $N$ by Theorem \ref{ub loop}. 
Then $\st_H(u)-\{f\}=\st_{H'}(u)$ for any $v\neq u\in V(H)$, implying that  $f\in \loops_{H'}(v)$ and $H-\{v,f\}=H'-\{v,\loops(H')\}$.

Assume that $H-\{v,f\}$ has an unbalanced cycle $C$. Then $E(C)\cup\{f\}\in\cC(N)$ as $H$ is  a lifted-graphic representation for $N$. On the other hand, since $C$ is a cycle of $H-\{f,v\}$ of length at least 2, $C$ is also an unbalanced cycle of $H'-\{f,v\}$ by \ref{H-v}. Since $H'$ is a frame representation for $N$ and $f\in \loops_{H'}(v)$, we have $E(C)\cup\{f\}\in\cI(N)$, a contradiction. Hence, $H-\{f,v\}$ is balanced. That is, the claim holds.
\end{proof}

By \ref{1 case} and symmetry, we may therefore assume that neither $H$ nor $H'$ is a lifted-graphic representation of $N$ with a loop. Then $\st_H(u)=\st_{H'}(u)$ for any $v\neq u\in V(H)$. Since $H\neq H'$, there is a link $e=vu$ of $H$ (or $H'$), which is a loop of $H'$ (or $H$) incident with $u$. Assume that $H-\{v,\loops(H)\}$ has an unbalanced cycle $C$. Since $|E(C)|\geq2$, it follows from \ref{H-v} that $C$ is also an unbalanced cycle of $H'-v$. Let $P$ be a minimal path in $H-v$ joining $u$ and $C$. Note that $P=\{u\}$ when $u\in V_H(C)$. Comparing $H[E(C\cup P)\cup\{e\}]$ and $H'[E(C\cup P)\cup\{e\}]$, we will get a contradiction. Hence, $H-\{v,\loops(H)\}$ is balanced. That is, (2) holds.
\end{proof}

In Section 5 we will need a number of simple conditions which prevent a matroid from being an excluded minor for the class of quasi-graphic matroids. In the following Lemmas we gather a few such conditions. Lemmas \ref{almost fix}-\ref{extendable} will be only used in the proof of Theorem \ref{main-result+1}.

\begin{lemma}\label{2-way-extendable}
Let $e,f$ be elements of a $3$-connected matroid $N$ such that $N\del e$, $N\del f$, and $N\del e,f$ are $3$-connected. Let $H$ be a $2$-connected unbalanced framework for $N\del e,f$ that has no blocking vertices. If $H$ can be extended to frameworks for $N\del e$ and $N\del f$, then $N$ is quasi-graphic.
\end{lemma}
\begin{proof}
Let $G$ be a graph with $H=G\del e,f$ such that $G\del e$ and $G\del f$ are frameworks for $N\del e$ and $N\del f$, respectively. Since $H$ is connected, by Lemmas \ref{3-con case} and \ref{2-conn} we may assume that $G$ is 2-connected. We claim that $G$ is a framework for $N$. Evidently, (QG1) and (QG2) hold. Since $N$ is 3-connected, by Theorem \ref{3-con case-1}, it suffices to show that (QG3) holds. Let $v$ be a vertex of $G$. When $e,f\in\st_{G}(v)$, (QG3) obviously holds for $v$. So by symmetry we may assume that $e\notin\st_G(v)$. Since $H$ is 2-connected and has no blocking vertices, $H-v$ is connected and unbalanced. Then it follows from Lemma \ref{e in C} that $e\in \cl_N(E(H-v))$ as $G\del f$ is a framework for $N\del f$. When $f\in\st_{G}(v)$, since $\cl_N(E(G-v))=\cl_N(E(H-v))$ and $G\del e$ is a framework for $N\del e$, (QG3) holds for $v$.
When  $f\notin\st_{G}(v)$, by the symmetry between $e$ and $f$, we have $f\in \cl_N(E(H-v))$. So $\cl_N(E(G-v))=\cl_N(E(H-v))$, implying that (QG3) holds for $v$ as  $\st_{G}(v)=\st_H(v)$.
\end{proof}

\begin{lemma}\label{extendable}
Let $e,f$ be elements of a $3$-connected matroid $N$ such that $N\del e$, $N\del f$, and $N\del e,f$ are $3$-connected. Let $H$ be a $2$-connected framework for $N\del e$. Assume that there is a balanced cycle $C$  of $H$ with $f\in E(C)$ such that all vertices in $V_H(C)$ are fixed in $H\del f$. If $N\del f$ is quasi-graphic, so is $N$.
\end{lemma}
\begin{proof}
Let $G^{''}$ be a framework for $N\del e,f$ that can be extended to a framework $G'$ for $N\del f$. By Lemmas \ref{3-con case} and \ref{2-conn} we may further assume that $G^{''}$ and $G'$ are 2-connected.
Since $G^{''}$ and $H\del f$ are frameworks for $N\del e,f$ and  all vertices in $V_H(C)$ are fixed in $H\del f$, by Lemma \ref{fix}, we may assume that corresponding vertices in $G^{''}[E(C)-\{f\}]$ and $H[E(C)-\{f\}]$ are labelled by same symbols and
\[\st_{G''}^*(v)=\st_{H\del f}^*(v)=\st^*_H(v)-\{f\}\in\cC^*(N\del e,f),\] for any $v\in V_H(C)$.
Hence, $G^{''}[E(C)-\{f\}]$ and $H[E(C)-\{f\}]$ are isomorphic paths. Let $G$ be the graph obtained from $G'$ by adding $f$ to $G'$ such that $G[E(C)]$ is a cycle. That is, $G[E(C)]$ and $H[E(C)]$ are isomorphic.

We claim that $G$ is a framework for $N$. (QG1) obviously holds. Since $G$ is 2-connected and $r(N)=r(N\del f)=|V(G')|$, (QG2) holds for $G$. Since $N$ is 3-connected, by Theorem \ref{3-con case-1} it suffices to show that (QG3) holds. Since $E(C)$ is a circuit of $N$ and $G'$ is a framework for $N\del f$, (QG3) holds  for each vertex in $V(G)-V_G(E(C))+V_G(f)$. For any $v\in V_G(E(C))-V_G(f)$, since $v$ is fixed in both $H$ and $H\del f$ by Lemma \ref{fix-0}, we have $$\st_{G'}^*(v)-\{e\}=\st_{G''}^*(v)=\st^*_H(v)\in\cC^*(N\del e,f)\cap\cC^*(N\del e) \eqno{(2.1)}$$ by Lemma \ref{fix}. Since $N\del(\st_{G''}^*(v)\cup\{e,f\})$ is 3-connected and non-graphic, $G'-v$ is unbalanced, so $\st_{G'}^*(v)\in\cC^*(N\del f)$. Combined with (2.1), $\st_{G'}^*(v)\in\cC^*(N)$ or $\{e,f\}\in\cC^*(N)$. Since $N$ is 3-connected, $\st_{G'}^*(v)\in\cC^*(N)$. Hence, (QG3) holds for $v$.
\end{proof}

\section{9-connected excluded minors with rank less than nine.}
In this section, we prove that, if $M$ is a 9-connected excluded minor for the class of quasi-graphic matroids with $r(M)\leq8$, then $M$ is isomorphic to $U_{3,7}$ or  $U_{4,7}$. To prove this, we need one more definition.

Let $G$ be a simple graph. For a positive integer $k$, let $kG$ denote the graph obtained from $G$ by replacing each edge of $G$ by a parallel class with exactly $k$ edges.

\begin{theorem}\label{U3,7}
$U_{3,7}$ is an excluded minor for the class of quasi-graphic matroids.
\end{theorem}
\begin{proof}
First we show that $2K_3$ is the unique framework for $U_{3,6}$.
Let $G$ be a framework for $U_{3,6}$. By Lemma \ref{3-con case}
we may assume that $G$ is connected. Then $|V(G)|=3$. Since $|E(G)|=6$,
either each vertex in $G$ is incident with exactly four edges
or some vertex $v$ is incident with at most three edges. When the former case happens,
$G$ is isomorphic to $2K_3$. When the latter case happens,
since $G-v$ has no balanced cycles with at most two edges,
Lemma \ref{circuit in framework} implies that $U_{3,6}$ has a triangle, a contradiction.

Since $2K_3$ is the unique framework for $U_{3,6}$, it is easy to verify that $U_{3,7}$ is not quasi-graphic. Moreover, since $6K_2$ is a framework for $U_{2,6}$, the theorem holds.
\end{proof}

\begin{theorem}\label{U4,7}
$U_{4,7}$ is an excluded minor for the class of quasi-graphic matroids.
\end{theorem}
\begin{proof}
Let $C_4$ be a 4-edge cycle, let $K$ be the graph obtained from $2C_4$ by deleting a pair of  non-adjacent edges. Evidently, $K_4$ and $K$ are frame representations for $U_{4,6}$. Note that, neither $K_4$ nor $K$ can be extended to a framework for $U_{4,7}$. Since $2K_3$ is a framework for $U_{3,6}$, to prove the theorem, it suffices to show that, besides $K_4$ and $K$, $U_{4,6}$ has no other frameworks.

Let $G$ be a framework for $U_{4,6}$. By Lemma \ref{3-con case} we may assume that $G$ is connected. Then $|V(G)|=4$, $G$ is 2-connected and each vertex of $G$ is incident with at least three edges. Assume that $G$ has a blocking vertex $u$. Since each circuit in $U_{4,6}$ has five elements and $G-u$ is balanced, $G-u$ is a forest, so $|\st_G(u)|\geq4$. Since $G$ is 2-connected, $G-u$ is a 2-edge path; that is, $|\st_G(u)|=4$. Let $v_1,v_2$ be the degree-1 vertices of $G-u$. Since $E(G)-\{f\}$ is a circuit of $U_{4,6}$ for each edge $f\in\st_G(u)$, there are exactly two edges joining $u$ and $v_i$ for each $1\leq i\leq 2$. So $\st_G(u)$ is dependent in $U_{4,6}$, a contradiction. So $G$ has no blocking vertices. For each vertex $v$ of $G$, since $G-v$ is connected and unbalanced, $|E(G-v)|=3$ as $|\st_G(v)|\geq3$. So $|\st_G(v)|=3$. Since $|E(G)|=6$, by the arbitrary choice of $v$, the graph $G$ has no loops and $G$ is isomorphic to $K_4$ or $K$.
\end{proof}

\begin{theorem}\label{rank<8}
Let $M$ be an excluded minor for the class of quasi-graphic matroids. If $M$ is $9$-connected with rank at most eight, then $M$ is isomorphic to $U_{3,7}$ or $U_{4,7}$.
\end{theorem}
\begin{proof}
We claim that $M$ is isomorphic to $U_{r,2r-1}$, $U_{r,2r}$, $U_{r,2r+1}$, or $U_{8,n}$ for a number $r$, where $n\geq 15$. Assume that $M$ has a circuit $C$ with $|C|\leq r(M)$. Without loss of generality we may further assume that $C$ is chosen as small as possible. When $|E(M)-C|\geq|C|$, the partition $(C,E(M)-C)$ is a $|C|$-separation, a contradiction to the fact that $M$ is 9-connected. When $|E(M)-C|<|C|$, since $|C|\leq r(M)$ and $E(M)-C$ is independent by the choice of $C$, the partition $(C,E(M)-C)$ is an $|E(M)-C|$-separation, a contradiction. So $M$ is uniform. Then it follows from (\cite{Oxley}, Corollary 8.6.3) that the claim holds.


Since $kK_2$ is a framework for $U_{2,k}$, we have $r(M)\geq3$. Since $U_{4,7}$ is a minor of $U_{r,2r-1}$, $U_{r,2r}$, $U_{r,2r+1}$, and $U_{8,n}$ when $r\geq4$ and $n\geq 15$, by Theorem \ref{U4,7} either $r(M)=3$ or $M$ is isomorphic to $U_{4,7}$. Moreover, since $U_{3,6}$ is quasi-graphic, the theorem holds from  Theorem \ref{U3,7} .
\end{proof}

\section{Frameworks for graphic matroids}
Let $G$ be a graph, and $M(G)$ its cycle matroid. A {\sl signed graph} is a pair $(G,\Sigma)$ with  $\Sigma\subseteq E(G)$, each edge in $\Sigma$ is labelled by $-1$ and other edges are labelled by 1. A cycle $C$ of $G$ is {\sl $\Sigma$-even} if $|E(C)\cap\Sigma|$ is even, otherwise it is {\sl $\Sigma$-odd}. A set $\Sigma'\subseteq E(G)$ is a {\sl signature} of $(G,\Sigma)$ if $(G,\Sigma)$ and $(G,\Sigma')$ have the same $\Sigma$-even cycles and the same $\Sigma$-odd cycles. Evidently, for any cut $C^*$ of $G$, the set $\Sigma\triangle E(C^*)$ is a signature of $(G,\Sigma)$. For a framework $H$ for a matroid $N$, we say that $H$ is a \emph{signed graph} if there is a set $\Sigma\subseteq E(H)$ such that a cycle $C$ of $H$ is balanced if and only if $C$ is a $\Sigma$-even cycle. We also say that $\Sigma$ is a {\sl signature} of $H$.

All definitions in the following five paragraphs were first given by Chen, DeVos, Funk and Pivotto \cite{CDFP}.

\smallskip
\emph{Fat thetas}. Let $G_1, G_2, G_3$ be non-empty graphs with distinct vertices $x_i, y_i \in
V(G_i)$. Let $G$ be obtained from $G_1, G_2, G_3$ by identifying
$y_i$ and $x_{i+1}$  to a vertex  $w_i$ for every $1\leq i\leq 3$ (where the indices are modulo $3$).
Let $H$ be obtained from $G_1,G_2,G_3$ by identifying $x_1, x_2, x_3$ to a
vertex $x$ and identifying $y_1, y_2, y_3$ to a vertex $y$. A cycle
of $H$ is balanced if and only if $E(C)$ is completely
contained in one of $G_1, G_2$ or $G_3$.
Then we say that $H$ is a \emph{fat theta} obtained from $G$.


\smallskip
\emph{Simple curlings}. Let $G$ be a graph and $v\in V(G)$. Let $H$ be the signed graph obtained from $G$ by first labelling all edges incident with $v$ by $-1$, and then changing any such edge $e=vu$ to  a loop incident with $u$ while keeping all other edges not incident with $v$ unchanged and labelled by $1$. Then we say that $H$ is a {\em simple curling} of $G$.


\smallskip
\emph{Pinches}. \ If $H$ is obtained from a graph $G$ by
identifying two vertices $v_1$ and $v_2$ to  a new vertex $v$ and
labeling all edges originally incident with $v_1$ by $-1$ and all other edges
by $1$, then we say $H$ is a \emph{pinch}. An edge with
ends $v_1, v_2$ becomes an unbalanced loop incident to $v$

\smallskip
\emph{$4$-twistings}. \ Let $G_1,G_2,G_3,G_4$ be graphs (not
necessarily all non-empty) with distinct vertices $x_i,y_i,z_i\in
V(G_i)$. Let $G$ be obtained from $G_1,G_2,G_3,G_4$ by identifying
$x_i,y_{3-i},z_{i+2}$ to a vertex $w_i$ for every $1\leq i\leq 4$ (where
the indices are modulo 4). Let $H$ be a signed graph obtained
from $G_1,G_2,G_3,G_4$ by identifying $x_1$, $x_2,x_3,x_4$ to a
vertex $x$, identifying $y_1,y_2,y_3,y_4$ to a vertex $y$ and
identifying $z_1,z_2,z_3,z_4$ to a vertex $z$, and with all edges
originally incident with $x_1, y_2$ or $z_3$ labelled by $-1$ and
all other edges labelled by $1$. 
Then we say that $H$ is a \emph{$4$-twisting} of $G$.

\smallskip
\emph{Consecutive twistings}. \
Let $G_1, \ldots, G_k$
(for $k \geq 3$), be graphs with distinct vertices $x_i, y_i, z_i
\in  V(G_i)$ for $1\leq i \leq k$. Let $G$ be a graph obtained from $G_1,
\ldots, G_k$ by identifying $z_1, z_2, \ldots, z_k$ to a vertex $z$
and for each $1\leq i \leq k$ identifying $y_{i-1}$ and $x_i$ to a vertex
$w_i$ (where the indices are modulo $k$). Let $H$ be the signed
graph obtained from $G_1, \ldots, G_k$ by identifying
$y_{i-1},z_i,x_{i+1}$ to a vertex $u_i$ for every $1\leq i\leq k$ (where
the indices are modulo $k$), and with all edges originally incident
with $y_1$ or $x_2$ labelled by $-1$ and all other edges labelled by $1$. 
Then we say that $H$ is a \emph{consecutive twisting} or a \emph{consecutive $k$-twisting} of $G$.
If $k$ is odd then $H$ is a \emph{consecutive odd-twisting} of $G$.

\begin{theorem}(\cite{CDFP}, Corollary 1.3.)\label{4-con}
Let $G$ be a $3$-connected graph with $|V(G)|\geq5$. Let $H$ be a frame representation of $M(G)$. Then either $H$ is balanced, or $H$ is obtained from $G$ as a simple curling, a pinch, a $4$-twisting, or a consecutive odd-twisting.
\end{theorem}

Recall that $c(H)$ is the number of components of $H$.

\begin{theorem}(\cite{GGW}, Theorem 2.7.)\label{determin bal}
Let $H$ be a framework for a matroid $N$. If $r(N)\leq |V(H)|-c(H)$, then $N=M(H)$.
\end{theorem}

\begin{theorem}(\cite{Zaslavsky87}, Theorem 2.)\label{Zas}
Let $H$ be a lifted-graphic representation of a matroid $N$. Then $N$ is binary if and only if $H$ is a signed graph or $H$ has a unique unbalanced component which is a fat theta.
\end{theorem}

Let $G$ be a graph, and let $(X_1, X_2)$ be a partition of $E(G)$
such that $V(X_1) \cap V(X_2) = \{u_1, u_2\}$. We say that $G'$
is obtained by a \emph{Whitney flip} of $G$ on $\{u_1,u_2\}$ if $G'$
is a graph obtained by identifying vertices $u_1,u_2$ of $G[X_1]$
with vertices $u_2,u_1$ of $G[X_2]$, respectively.
A graph $G'$ is \emph{$2$-isomorphic to} $G$ if $G'$ is obtained from
$G$ by a sequence of the operations: Whitney flips, identifying two
vertices from distinct components of a graph, or partitioning a
graph into components each of which is a block of the original
graph.

In his Ph.D. thesis, Shih \cite{Shih} proved the following characterization of graphic lifted-graphic matroids (see also \cite{Pivotto},
Theorem 4.1.).

\begin{theorem}[Theorem 1, Chapter 2 in \cite{Shih}]\label{Shih-Pivotto}
Let $G$ be a graph and let $H$ be a lifted-graphic representation of $M(G)$. Assume that $H$ is an unbalanced signed graph. Then there exists a graph $G'$ $2$-isomorphic to $G$ such that one of the following holds.
\begin{itemize}
    \item[$(1)$]  $H$ is obtained from $G'$ by a pinch.
    \item[$(2)$]  $H$ is obtained from $G'$ by a $4$-twisting.
    \item[$(3)$] $H$ is obtained from $G'$ by a consecutive twisting.
\end{itemize}
\end{theorem}

Following a similar way as the proof of  (\cite{GGW}, Theorem 1.4.), we prove

\begin{theorem}\label{repre matroid}
Let $H$ be a $2$-connected framework for a $3$-connected matroid $N$. If $N$ is representable, then $H$ is a frame representation or a lifted-graphic representation of $N$.
\end{theorem}
\begin{proof}
Without loss of generality we may assume that $H$ is unbalanced. Then $|V(H)|=r(N)$. By Theorem \ref{ub loop} we may assume that $H$ has no loops. Assume that there is a vertex $v$ of $H$ such that $r_N(E(H-v))\leq r(N)-2$. Since $H-v$ is connected, it follows from Theorem \ref{determin bal} that $H-v$ is balanced. Then $v$ is a blocking vertex of $H$, so $H$ is a frame representation and a lifted-graphic representation of $N$. So we may assume that $r_N(E(H-v))=r(N)-1$ for each vertex $v$ of $H$. Moreover, since $H$ has no loops, $\st_H(v)$ is a cocircuit of $N$ by (QG3).

Let $A$ be a matrix over a field $\mathbb{F}$ with linearly independent rows satisfying $N=M(A)$, where $M(A)$ is the matroid represented by $A$. Since $\st_H(v)$ is a cocircuit of $N$ for each vertex $v$ of $H$, there is a matrix $B\in\mathbb{F}^{V(H)\times E(H)}$ such that
\begin{enumerate}
\item the row-space of $B$ is contained in the row-space of $A$, and
\item for each $v\in V(H)$ and $e\in E(H)$, the element of $B$ in the row labelled by $v$ and the column labelled by $e$ is non-zero if and only if $v$ is incident with $e$ in $H$.
\end{enumerate}
Note that $M(B)$ is a frame matroid and $H$ is a framework for $M(B)$. Since $H$ is connected, we have that $$|V(H)|=r(M(A))\geq r(M(B))\geq |V(H)|-1, \eqno{(4.1)}$$ and if $r(M(A))=r(M(B))$ then $M(A)=M(B)$ by (1) and (2). So we may assume that $r(M(A))>r(M(B))$. Then $M(B)=M(H)$ by Theorem \ref{determin bal}, up to row-operations we may assume that $A$ is obtained from $B$ by appending a single row by (4.1). Hence, $H$ is a lifted-graphic representation of $N$.
\end{proof}

By Lemma \ref{3-con case} (or Lemma \ref{disconnected}) and Theorem \ref{repre matroid}, we have

\begin{corollary}\label{repre matroid+1}
Let $H$ be a framework for a $3$-connected representable matroid $N$. Then $H$ is a frame representation or a lifted-graphic representation of $N$.
\end{corollary}

The following result is an immediate consequence of Theorems \ref{4-con}, \ref{Zas}, \ref{Shih-Pivotto}, and Corollary \ref{repre matroid+1}.

\begin{theorem}\label{graphic matroid}
Let $G$ be a $3$-connected graph with $|V(G)|\geq5$, and $H$ a connected framework for $M(G)$. Then $H$ is isomorphic to $G$, or $H$ is obtained from $G$ by a simple curling, a pinch, a $4$-twisting, or a consecutive twisting.
\end{theorem}

By Theorem \ref{graphic matroid} we have

\begin{corollary}\label{4-con graphic matroid}
Let $G$ be a $4$-connected graph with $|V(G)|\geq5$, and $H$ a connected framework for $M(G)$. Then $H$ is isomorphic to $G$ or $H$ is obtained from $G$ by a simple curling or a pinch.
\end{corollary}

\begin{lemma}\label{graphic matroid+1}
Let $G$ be a $3$-connected simple graph, and $H$ a $4$-connected unbalanced framework for $M(G)$ with $|V(H)|\geq4$. Then
\begin{enumerate}
\item $H$ is obtained from $G$ by a simple curling or a pinch, or
\item $H$ is a signed graph with a signature $X$ such that $H[X]$ is a triangle. 
\end{enumerate}
\end{lemma}
\begin{proof}
Assume that (1) is not true. Since $H$ is a 4-connected unbalanced graph with $|V(H)|\geq4$, by Theorem \ref{graphic matroid}, the graph $H$ is obtained from $G$ by a 4-twisting or a consecutive 3-twisting. Without loss of generality that it is a 4-twisting, since the consecutive 3-twisting is (up to relabelling of vertices) the special case of this in which one of the $G_i$ has no edges. By symmetry we may assume without loss of generality that none of $G_1, G_2, G_3$ has more than 3 vertices, where $G_i$ and symbols that will be used but not defined in the proof, say $w_i,x_i, y_i, x, y, z$, are defined as in the definition of 4-twistings. By 3-connectivity and simplicity of $G$, there is precisely one edge $e_i$ from $w_1$ to each $w_i$ with $i\in \{2,3,4\}$, and there are no other edges incident with $w_1$. By the definition of 4-twistings, the signature of $H$ is $\{e_2, e_3, e_4\}$. The edge $e_2$ can only arise from an edge $x_1y_1$ in $G_1$ or an edge $x_2y_2$ in $G_2$: in either case it joins $x$ to $y$ in $H$. Similarly $e_3$ joins $x$ to $z$ in $H$ and $e_4$ joins $y$ to $z$ in $H$. Thus the signature of $H$ is the set of edges of a triangle.
\end{proof}


\section{Proof of Theorem \ref{main-result}.}
Recall that $c(H)$ is the number of components of a graph $H$. Lemmas \ref{2.1}-\ref{non-3-con} will be frequently used in this section.

\begin{lemma}\label{2.1}
Let $H$ be a framework for a matroid. For an edge $f\in E(H)$, if $H\del f$ is balanced and $H$ has a balanced cycle containing $f$, then $H$ is balanced.
\end{lemma}
\begin{proof}
Since $r(E(H))=r(E(H\del f))=|V(H)|-c(H)$, the graph $H$ is balanced.
\end{proof}

Note that Lemma \ref{2.1} also follows immediately from the theta property.

For any subset $X$ of $E(H)\cup V(H)$, if $H\del X$ is balanced, we say that $X$ is a {\sl balancing set}  of $H$. 
Note that, when $H$ is balanced, each subset of $E(H)\cup V(H)$ is balancing. We say a balancing set $X$ is {\em minimal} if no proper subset of $X$ is a balancing set of $H$. Note that, when $H$ has a non-empty minimal balancing set $V\cup E$ with $V\subseteq V(H)$ and $E\subseteq E(H)$, the graph $H$ is unbalanced and $E\cap \st(V)=\emptyset$ by the definition of minimal balancing sets.

\begin{lemma}\label{minimal balancing}
Let $H$ be a connected unbalanced framework for a matroid $N$. If $X$ is a minimal balancing edge set of $H$ with $X\subseteq E(H)$, then $X\in\cC^*(N)$.
\end{lemma}
\begin{proof}
Since $H\del X$ is connected and balanced, $r(E(H\del X))=|V(H)|-1=r(N)-1$. On the other hand, since each cycle in $H\del (X\del\{f\})$ containing $f$ is unbalanced for each $f\in X$ by Lemma \ref{2.1},  $r(E(H\del X)\cup\{f\})=r(N)$. Hence, $X$ is a cocircuit of $N$.
\end{proof}

\begin{lemma}\label{rank of balance set}
Let $H$  be a connected framework for an $n$-connected matroid $N$ with $|V(H)|\geq n$.  When $H$ is unbalanced, each balancing set of $H$ that contains only edges has rank at least $n$.
\end{lemma}
\begin{proof}
Assume not. Let $X$ be a minimal balancing set of $H$ with $X\subseteq E(H)$ and $r(X)\leq n-1$. Then $r(E(N)-X)=r(N)-1$ by Lemma \ref{minimal balancing}. Since $H\del X$ contains a spanning tree of $H$, $|E(N)-X|\geq n-1$, so $(X,E(N)-X)$ is an $r(X)$-separation of $N$, a contradiction.
\end{proof}


For any subset $X$ of $E(H)\cup V(H)$, if $c(H\del X)>c(H)$, then we say that $X$ is a {\em cut} of $H$.

\begin{lemma}\label{5.2}
Let $H$ be a framework for a matroid. Let $X_i=V_i\cup E_i$ be a balancing set of $H$ with $V_i\subseteq V(H)$ and $E_i\subseteq E(H)$ for each $1\leq i\leq 2$.
\begin{enumerate}
\item If $X_1$ is minimal and contains a link $f$ satisfying $f\notin E_2\cup\st(V_2)$, then $X_1\cup X_2$ contains a cut of $H$.
\item If $H-(X_1\cup X_2)$ is connected and $V_1\cap V_2=\emptyset$, then $E_1\cup E_2\cup E(H[V_1\cup V_2])$ is a balancing set of $H$.
\end{enumerate}

\end{lemma}
\begin{proof}
First we prove that (1) is true.  Since $X_1$ is minimal, each cycle in $H\del (X_1\del\{f\})$ containing $f$ is unbalanced by Lemma \ref{2.1}. Moreover, since $H-X_2$ is balanced and $f\in E(H-X_2)$, the graph $H\del (X_1\del\{f\})$ has a cut contained in $X_2\cup\{f\}$, so (1) holds.

Assume that (2) is not true. Let $C$ be an unbalanced cycle of $H\del(E_1\cup E_2\cup E(H[V_1\cup V_2]))$ with $|V(C)\cap (V_1\cup V_2)|$ as small as possible. Since $H-X_i$ is balanced for each $1\leq i\leq2$, we have $V(C)\cap V_i\neq\emptyset$. Then $|V(C)\cap (V_1\cup V_2)|\geq2$ as $V_1\cap V_2=\emptyset$. Since $C$ does not contain edges in $ H[V_1\cup V_2]$, the subgraph $C-(V_1\cup V_2)$ is disconnected.  Moreover, since $H-(X_1\cup X_2)$ is connected, there is a path $P$ of $H-(X_1\cup X_2)$ connecting two components of $C-(V_1\cup V_2)$ such that $C\cup P$ is a theta subgraph. For any cycle $C'$ of $H$ with $P\subseteq C'\subseteq C\cup P$, since $|V(C')\cap (V_1\cup V_2)|\leq|V(C)\cap (V_1\cup V_2)|-1$, the cycle $C'$ is balanced by the choice of $C$. Therefore $C$ is balanced by the theta property, a contradiction. So (2) holds.
\end{proof}

Note that, the set $E_i$ in Lemma \ref{5.2} may be empty.

Let $X$ and $Y$ be subsets of the ground set of a matroid $N$. Set  \[\sqcap_N(X,Y)=r_N(X)+r_N(Y)-r_N(X\cup Y).\] When $(X,Y)$ is a partition of $E(N)$, we often denote $\sqcap_N(X,Y)$ by $\lambda_N(X)$. When there is no confusion, subscripts will be omitted.

\begin{lemma}\label{non-3-con}
Let $H$ be a $4$-connected framework for a simple and non-$3$-connected matroid $N$ with $|V(H)|\geq4$. Then $H$ is unbalanced and has a balancing set $X$ with $r(X)\leq 2$. In particular, when $N$ has no triangles, $|X|\leq 2$.
\end{lemma}
\begin{proof}
Since $N$ is not 3-connected and $H$ is 4-connected, $H$ is unbalanced. Let$(X,Y)$ be an exact $k$-separation of $N$ for an integer $1\leq k\leq2$. We may assume that $(X,Y)$ is chosen with $\lambda(X)+c(H[X])+c(H[Y])$ as small as possible.

\vspace{3mm}
Case 1. $H[X]$ and $H[Y]$ are connected.
\vspace{3mm}

Set $m=|V_H(X)\cap V_H(Y)|$. Then $m\in\{k-1,k,k+1\}$ as $\lambda(X)=k-1$. Since $H$ is 4-connected, by symmetry we may assume that $V_H(Y)=V(H)$ and $m=|V_H(X)|$. When $m=k-1$, we have that $k=2$ and  $H[X]$ consists of loops, so $N$ has a circuit contained in $X$ of size at most 2, which is not possible as $N$ is simple. When $m=k+1\leq3$, both $H[X]$ and $H[Y]$ are balanced, that is, $X$ is a balancing set of $H$ with $r(X)=k$. When $m=k\leq2$,  one of $H[X]$ and $H[Y]$ is balanced and the other is unbalanced. If $H[X]$ is balanced, since $|X|\geq k$ and $k=|V_H(X)|$, the set $X$ contains a circuit of $N$ whose size  is at most 2, a contradiction. So $H[Y]$ is balanced. That is, $X$ is a balancing set of $H$ with $r(X)=k-1$.

\vspace{3mm}
Case 2. $H[X]$ is disconnected, implying $|X|\geq2$.
\vspace{3mm}

Let $X_1$ be the edge set of a component of $H[X]$.
\begin{claim}\label{local conn}
Either $|X-X_1|< k$ or $\sqcap(X_1, X-X_1)=1$ and $\sqcap(X_1,Y)=0$.
\end{claim}
\begin{proof}[Subproof.]
Assume that $|X-X_1|\geq k$. Since $c(H[X])+c(H[Y])> c(H[X-X_1])+c(H[Y\cup X_1])$, we have $\lambda(X-X_1)>\lambda(X)$ by the choice of $(X,Y)$, so $\sqcap(X_1, X-X_1)>\sqcap(X_1,Y)\geq0$. Since $\sqcap(X_1, X-X_1)\leq1$, the claim holds.
\end{proof}

Assume that $\sqcap(X_1, X-X_1)=0$. By \ref{local conn}, we have $1\leq |X-X_1|<k$. So $k=2=c(H[X])$. Using \ref{local conn} again, we have $|X_1|=1$, so $|X|=2$, implying that $H[Y]$ is a connected spanning subgraph of $H$. Since $N$ is simple, $r(X)=2$. Then $H[Y]$ is balanced as $\lambda(X)=1$, so the lemma holds. Hence, we may assume that $\sqcap(X', X-X')=1$ for the edge set $X'$ of each component of $H[X]$, implying that $H[X']$ is unbalanced by Lemma \ref{circuit in framework}. By symmetry we may further assume that (a) either $H[Y]$ is connected or each component of $H[Y]$ is unbalanced.

When $|X|=2$, since $\sqcap(X_1, X-X_1)=1$, we have $X\in\cC(N)$, which is not possible as $N$ is simple. So $|X|\geq3$. By symmetry assume that $|X-X_1|\geq2$ . By \ref{local conn}, we have $\sqcap(X_1,Y)=0$. Since $X_1$ is unbalanced, each component of $H[Y]$ that shares vertices with $H[X_1]$ is balanced by Lemma \ref{circuit in framework}. Then $H[Y]$ is connected and balanced by (a), implying that $r(Y)=|V_H(Y)|-1$ and $X$ is a balancing set of $H$.

Let $X_1, \ldots, X_{c(H[X])}$ be the components of $H[X]$. Since an unbalanced spanning unicyclic subgraph of $H[X]$ is an independent set in $N$, we have $r(X)\geq\Sigma_1^{c(H[X])} |V_H(X_i)|-c(H[X])+1.$ Then $$1\geq k-1=\lambda(X)\geq |V_H(X)\cap V_H(Y)|-c(H[X])$$ as $r(Y)=|V_H(Y)|-1$. Hence,
$|V_H(X_i)\cap V_H(Y)|\leq2$ for each $1\leq i\leq c(H[X])$ and at most one $|V_H(X_i)\cap V_H(Y)|$ is not equal to 1. Since $H$ is 4-connected and $H[X]$ is disconnected, $V_H(Y)=V(H)$. Since $r(Y)=r(N)-1$ and  $\lambda(X)\leq1$, we have $r(X)\leq2$. Hence, the lemma holds as $X$ is a balancing set of $H$.
\end{proof}

Recall that we define $H'=H-\loops(H)$ when $H$ is a lifted-graphic representation of $N$, otherwise let $H'=H$. For $v_1,v_2\in V(H)$, we say $\{v_1,v_2\}$ is a {\em blocking pair} of $H$ if $v_i$ is a blocking vertex of $H'-v_{3-i}$ for each $1\leq i\leq 2$. Note that, by our definition, balanced frameworks have no blocking vertices, and no vertex in a blocking pair is a blocking vertex.

\begin{lemma}\label{case>2+}
Let $H$  be a $7$-connected unbalanced framework for a matroid with $|V(H)|\geq8$. Assume that $H$ has no blocking pairs and $H-\loops(H)$ has no blocking vertices. Then there is an edge $f$ of $H$ such that $H\del f$ has no blocking pairs and $H-(\loops(H)\cup\{f\})$ has no blocking vertices.
\end{lemma}
\begin{proof}
Let $e$ be an edge of $H$ and $S_e$ be a minimal subset of $V(H)$ such that $H-(\loops(H)\cup \{e\}\cup S_e)$ is balanced. We can further assume that $|S_e|\leq2$ otherwise the lemma holds. Let $f$ be a link of $H-(S_e\cup \{e\})$. Assume that $H\del f$ has a blocking pair $S_f$. Since $S_f\cup \{f\}$ is a minimal balancing set of $H$, it follows from Lemma \ref{5.2} (1) that $S_e\cup S_f\cup \{e,f\}$ contains a cut of $H$, a contradiction to the fact that $H$ is 7-connected, a contradiction. Following a similar way, we show that $H\del(\loops(H)\cup \{f\})$ has no blocking vertices. Hence, the lemma holds for $f$.
\end{proof}

\begin{lemma}\label{case>2}
Let $H$ be a $6$-connected framework for a $7$-connected matroid $N$ with $|V(H)|\geq7$. Assume that $H$ has no blocking pairs and $H-\loops(H)$ has no blocking vertices. Then at most one vertex of $H$ is not fixed.
\end{lemma}
\begin{proof}
Assume not. Let $v_1, v_2$ be unfixed vertices of $H$. Then $N\del \st^*_H(v_i)$ is graphic or non-3-connected for each $1\leq i\leq 2$. Since $H$ has no blocking pairs and $H-\loops(H)$ has no blocking vertices, by Lemma \ref{graphic matroid+1} or Lemma \ref{non-3-con}, for each $1\leq i\leq 2$, there is a minimal balancing set $X_i$ of $H\del\st^*_H(v_i)$ such that the following (a) or (b) happens. (a) $|X_i|=2$, the two edges in $X_i$ have no common vertex and if $X_i$ contains a loop then $H$ is not a lifted-graphic representation of $N$, for otherwise $H$ has a blocking pair by the definition of blocking pairs. (b) $H[X_i]$ is a triangle. Hence, $v_i\notin V(X_i)$ for each $1\leq i\leq 2$ no matter which case happens. Since $N$ is 7-connected, $X_i\cup\{v_i\}$ is a minimal balancing set of $H$ for each $1\leq i\leq 2$. 

We claim that $X_1=X_2$ when $v_2\notin V(X_1)$. When $X_1\subseteq X_2$, it follows from (a) and (b) that $X_1=X_2$. Hence, it suffices to show that $X_1\subseteq X_2$. Assume to the contrary that there is an edge $x\in X_1-X_2$. Then $x$ is not a loop of $H$, otherwise $x\in X_2$ as $v_2\notin V(X_1)$. Since $X_1\cup\{v_1\}$ is a minimal balancing set of $H$, by Lemma \ref{5.2} (1), $X_1\cup X_2\cup\{v_1,v_2\}$ contains a cut of $H$. Since $H$ is 6-connected, $X_1\cup X_2$ is a matching in $H$ of size 4 and $X_1\cup X_2\cup\{v_1,v_2\}$ is cut of $H$ by (a) and (b).
Let $H_1$ be a component of $H-\{v_1,v_2,X_1,X_2\}$. Set $H^+_1=H[V(H_1)\cup\{v_1,v_2\}]\del E(H[\{v_1,v_2\}])$. Since each balancing set of edges in $H$ has size at least 7 by Lemma \ref{rank of balance set}, $H\del (X_1\cup X_2)$ is unbalanced with $v_1, v_2$ as its blocking vertices. Then $H^+_1$ is balanced, otherwise $H^+_1$ has an unbalanced cycle containing exactly one vertex of $\{v_1,v_2\}$ by the theta property. So $\lambda(E(H^+_1))\leq 5$, a contradiction to the fact that $N$ is 7-connected.

When $X_1=X_2$, since $H-\{v_1, v_2, X_1\}$ is connected, by Lemma \ref{5.2} (2), $X_1\cup E(H[\{v_1,v_2\}])$ is a balancing set of rank at most 5, a contradiction to Lemma \ref{rank of balance set}. So $X_1\neq X_2$. By symmetry and the claim proved in the last paragraph, $v_i\in V(X_{3-i})$ for each $1\leq i\leq 2$. Let $x$ be the edge in $X_1$ that is not incident with $v_2$. Since each cycle in the 4-connected graph $H-\{v_1, v_2\}$ containing $x$ is unbalanced by Lemma \ref{2.1}, $x\in X_1\cap X_2$. Combined with (a) and (b), we have $r(X_1\cup X_2\cup E(H[\{v_1,v_2\}]))\leq6$. On the other hand, since $v_i\in V(X_{3-i})$, the graph $H-\{v_1, v_2, X_1, X_2\}$ is 2-connected. By Lemma \ref{5.2} (2) again, $X_1\cup X_2\cup E(H[\{v_1,v_2\}])$ is a balancing set, a contradiction to Lemma \ref{rank of balance set}.
\end{proof}

To prevent a matroid from being an excluded minor for the class of quasi-graphic matroids, we can use Lemmas \ref{case>2+} and \ref{case>2}  to show that $H$ and $H\del f$ have enough fixed vertices for some edge $f$, as long as $H-\loops(H)$ has no blocking vertex and $H$ has no blocking pair. In the rest of this section, we will show that when $H$ has a blocking vertex or a blocking pair, there is a balanced cycle $C$ of $H$ and $f\in E(C)$ such that all vertices in $V(C)$ are fixed in both $H\del f$ and $H$ (namely in Lemmas \ref{case:almost-balance} and \ref{case<3}). The case that $H$ has a blocking vertex will be dealt with first.

A biased graph $H$ is {\sl contra-balanced} if each cycle of $H$ is unbalanced.

\begin{lemma}(\cite{Zas81}, Theorem 6.)\label{signed-graph}
A biased graph is a signed graph if and only if it has no contra-balanced theta-subgraphs.
\end{lemma}

\begin{lemma}\label{signed-graph+1}
Let $v$ be a blocking vertex of a biased graph $H$. Let $x\in V(H-v)\cup E(H-v)$ that is not adjacent with $v$ when $x\in V(H-v)$. If $H-\{x,v\}$ is connected, then $H$ is a signed graph if and only if $H-x$ is a signed graph.
\end{lemma}
\begin{proof}
Evidently, it suffices to show that ``if" part. Assume that $H-x$ is a signed graph but $H$ is not. Then $H$ has a contra-balanced theta subgraph $T$ containing $x$ by Lemma \ref{signed-graph}. Since $v$ is a blocking vertex of $H$, $v$ is a degree-3 vertex of $T$. Since $v$ and $x$ are not adjacent when $x\in V(H-v)$, the graph $T-\{x,v\}$ has exactly two or three components. Since $H-\{x,v\}$ is connected, $H-\{x,v\}$ has a minimal forest $P$ that joins different components of $T-\{x,v\}$ such that $(T\cup P)-\{x,v\}$ is connected. Then $(T\cup P)-x$ consists of a theta subgraph $T'$ and some vertex-disjoint paths that are not in any cycle. Since $T$ is contra-balanced and $H-v$ is balanced, $T'$ is also a contra-balanced theta-subgraph by theta property. Hence, $H-x$ is not a signed graph by Lemma \ref{signed-graph}, a contradiction.
\end{proof}

Suppose that $v$ is a blocking vertex of $H$ and $H-v$ is connected.
In this case we define a relation $\sim_v$ on the edges in $\st_{H}(v)-\loops_H(v)$ by declaring $e \sim_v f$ if either $e=f$ or all cycles containing $e$ and $f$ are balanced. This is an equivalence relation, as we show next. Let $e_1,e_2,e_3$ be distinct edges in $\st_{H}(v)$ with $e_1 \sim_v e_2$ and $e_2 \sim_v e_3$. Let $T$ be a theta subgraph of $H$ containing all of $e_1, e_2$ and $e_3$; such a theta subgraph exists because $H-v$ is connected. The cycle in $T$ containing both $e_1$ and $e_2$ is balanced, and so is the cycle containing both $e_2$ and $e_3$. Therefore the cycle $C$ in $T$ containing $e_1$ and $e_3$ is balanced. Any other cycle containing $e_1$ and $e_3$ may be obtained from $C$ by rerouting along balanced cycles (contained in $H - v$), hence all the cycles containing $e_1$ and $e_3$ are balanced and $e_1 \sim_v e_3$, showing that $\sim_v$ is an equivalence relation. The same argument shows that a cycle of $H$ (that is not a loop) is unbalanced if and only if it contains two edges in $\st_{H}(v)$ which are not equivalent. We call the partition given by the equivalence classes of $\sim_v$ the {\em standard partition} of $\st_{H}(v)-\loops_H(v)$.  For more details, the reader can refer to (\cite{CP18}, Section 2) or (\cite{DF18}, Section 1). Definitions and results introduced in this paragraph will only be used in the proof of Lemma \ref{case:almost-balance-1}.

When $H$ is a signed graph with a blocking vertex $v$, since $H$ has no contra-balanced theta subgraph by Lemma \ref{signed-graph}, it is easy to show that 
there is a partition $(X_1,X_2)$ of $\st_H(v)-\loops_H(v)$ such that $H-(\loops(H)\cup X_i)$ is balanced for each $1\leq i\leq 2$. 
Split $v$ into $v_1, v_2$ such that $X_1, X_2$ are incident with $v_1$ and $v_2$, respectively, and such that each unbalanced loop in $H$ joins $v_1$ and $v_2$ and each balanced loop in $\loops_H(v)$ is a loop incident with any $v_i$. Let $G$ denote the new graph. Then $H$ is  a lifted-graphic representation of $M(G)$. Hence, if a framework for a matroid $N$ is a signed graph with a blocking vertex, then $N$ is graphic. This fact will be frequently used in the rest of this section without reference.

Let $X, Y\subseteq E(H)$ and $\mathcal{P}=(P_1,\ldots,P_n)$ be a partition of $X$. We will let $\mathcal{P}-Y$ denote the partition $(P_1-Y,\ldots,P_n-Y)$ of $X-Y$.

\begin{lemma}\label{case:almost-balance-1}
Let $H$ be a $5$-connected framework for a $5$-connected matroid $N$ with $|V(H)|\geq5$. Assume that $N$ is non-graphic and $H$ has a blocking vertex  $v$. Then the following hold.
\begin{enumerate}
\item For each vertex $v\neq u\in V(H)$, the graph $H-u$ is unbalanced and $N\del \st_{H}(u)$ is $3$-connected.
\item A vertex $u$ with $u\neq v$ is not fixed in $H$ if and only if $H\del E(H[\{v,u\}])$ is an unbalanced signed graph.
\item At most one vertex in $V(H)-\{v\}$ is not fixed in $H$.
\end{enumerate}
\end{lemma}
\begin{proof}
First we prove that (1) is true. If $H-u$ is balanced, then it follows from Lemma \ref{5.2} (2) that $E(H[\{u,v\}])$ is a balancing set of $H$ with rank at most 2, which is not possible by Lemmas \ref{minimal balancing} and \ref{rank of balance set}. Hence, $H-u$ is unbalanced.

Assume that $N\del \st_{H}(u)$ is not $3$-connected for some vertex $v\neq u\in V(H)$. Since $N$ has no triangles, by Lemma \ref{non-3-con}, $H-u$ has a minimal balancing set $X$  with $|X|\leq 2$. Since each cycle of $H-u$ containing exactly one edge of $X$ is unbalanced by Lemma \ref{2.1}, $X\subseteq \st_{H}(v)\cup\loops(H)$. Moreover, since $H-\{u,v\}$ is connected, it follows from Lemma \ref{5.2} (2) that $X\cup E(H[\{u,v\}])$ is a balancing set of $H$ with rank at most 4, which is not possible by Lemmas \ref{minimal balancing} and \ref{rank of balance set}.

Secondly, we prove that (2) is true. When $H\del E(H[\{v,u\}])$ is an unbalanced signed graph with $v$ as its blocking vertex, $N\del E(H[\{v,u\}])$ is graphic, and thus so is $N\del\st_H(u)$. Hence, $u$ is not fixed in $H$. Next, we prove the ``only if" part of (2) is true. Since $u$ is not fixed, $N\del \st_{H}(u)$ is graphic by (1). Then $H-u$ is an unbalanced signed graph with $v$ as its blocking vertex by (1), hence so is $H\del E(H[\{v,u\}])$ by repeatedly using Lemma \ref{signed-graph+1}.

Thirdly, we prove that (3) is true. Assume not. There are vertices $u,u'$ in $V(H)-\{v\}$ such that $H\del E(H[\{v,u\}])$ and $H\del E(H[\{v,u'\}])$ are signed graphs with $v$ as their blocking vertex by (2). Let $\mathcal{P}$ be the standard partition of $\st_{H}(v)-\loops_H(v)$ in $H$. Then $\mathcal{P}-\st_{H}(u)$ and $\mathcal{P}-\st_{H}(u')$ have exactly two non-empty members. When $\mathcal{P}-(\st_{H}(u)\cup\st_{H}(u'))$ has exactly two non-empty members, $\mathcal{P}$ has exactly two non-empty members, implying that $H$ is a signed-graph, a contradiction to the fact that $N$ is non-graphic. When $\mathcal{P}-(\st_{H}(u)\cup\st_{H}(u'))$ has exactly one non-empty member, since $N$ is simple, $H$ has a minimal balancing set $\{e,e'\}$ with $e\in E(H[\{v,u\}])$ and $e'\in E(H[\{v,u'\}])$. Then $\{e,e'\}\in\cC^*(N)$ by Lemma \ref{minimal balancing}, a contradiction.
\end{proof}

\begin{lemma}\label{case:almost-balance}
Let $H$ be a $6$-connected framework for a $6$-connected matroid $N$ with $|V(H)|\geq6$. Assume that $N$ is non-graphic and $H$ has a blocking vertex  $v$. Then $H-v$ has a balanced cycle $C$ such that all vertices in $V(C)$ are fixed in both $H\del f$ and $H$ for each edge $f$ of $C$.
\end{lemma}
\begin{proof}
We claim that a vertex $u\in V(H-v)$ is not fixed in $H$ if and only if it is not fixed in $H\del f$ for an arbitrary $f\in E(H-v)$. Evidently, it suffices to show that the ``if" part is true. Let $u$ be a vertex in $V(H-v)$ that is not fixed in $H\del f$. By Lemma \ref{case:almost-balance-1} (2), $H\del (E(H[\{v,u\}])\cup\{f\})$ is an unbalanced signed graph with $v$ as its blocking vertex. Since $H\del (E(H[\{v,u\}])\cup\{f,v\})$ is connected, $H\del E(H[\{v,u\}])$ is a signed graph by Lemma \ref{signed-graph+1}. So $u$ is not fixed in $H$ by Lemma \ref{case:almost-balance-1} (2) again.

By Lemma \ref{case:almost-balance-1} (3), $H-v$ has a balanced cycle $C$ such that all vertices in $V(C)$ are fixed in $H$. By the claim proved in the last paragraph, for each $f\in E(C)$, all vertices in $V(C)$ are also fixed in $H\del f$.
\end{proof}



Next, the case that $H$ has a blocking pair $S$ but $H-\loops(H)$ has no blocking vertices will be dealt with. To deal with this case, we need to introduce a characterization of the structure of biased graphs that have at least two blocking vertices. Lemma \ref{blocking vertex} will be only used in the proof of Lemma \ref{case<3-1}.

\begin{lemma}(\cite{Zas87}, Corollary 2.)\label{blocking vertex}
Let $V^*=\{v_1,\ldots,v_n\}$ be the set of blocking vertices of a biased graph $H$. Assume that $n\geq2$. Then one of the following holds.
\begin{enumerate}
\item $H$ is obtained from $mK_2$ by replacing each edge $e_i$ with a balanced graph $H_i$ such that all cycles of $H$ not contained in some $H_i$ are unbalanced, where $m\geq2$.
\item $H$ is obtained from a cycle $v_1v_2\ldots v_nv_1$ by replacing each edge $v_iv_{i+1}$ with a graph $H_i$ and a cycle in $H$ is unbalanced if and only if it contains $\{v_1,\ldots,v_n\}$, where no vertex in $H_i$ separates $v_i$ and $v_{i+1}$ and all subscripts are modulo $n$.
\end{enumerate}
\end{lemma}

\begin{lemma}\label{case<3-1}
Let $H$  be a $5$-connected unbalanced framework for a $6$-connected matroid $N$ with $|V(H)|\geq6$. Assume that $H$ has a blocking pair $S$ and $H-\loops(H)$ has no blocking vertices. Then we have
\begin{enumerate}
\item If some $v\in V(H)-S$ is not fixed in $H$, then there is a vertex $u\in S$ such that $\{u,v\}$ is a blocking pair of $H$.
\item $H$ has at most two blocking pairs, and they  have a common vertex.
\end{enumerate}
\end{lemma}

Note that at most one vertex in $V(H)-S$ can be contained in a blocking pair of $H$ and at most one vertex in $V(H)-S$ is not fixed in $H$ by Lemma \ref{case<3-1}.

\begin{proof}[Proof of Lemma \ref{case<3-1}.]
When $H$ is a lifted-graphic representation for $N$ with a loop, we may assume that the loop is in $\st(S)$. 

\begin{claim}\label{10.1}
If $S_1, S_2$ are blocking pairs of $H$, then $S_1\cap S_2\neq \emptyset$.
\end{claim}
\begin{proof}[Subproof.]
Assume otherwise. By Lemma \ref{5.2} (2), $E(H[S_1\cup S_2])$ is a balancing set of $H$ of rank at most 4, a contradiction to Lemma \ref{rank of balance set}.
\end{proof}

First, we prove that (1) is true. Since $N\del \st_H(v)$ is graphic or non-3-connected, by Lemma \ref{graphic matroid+1} or Lemma \ref{non-3-con}, either (a) $H-v$ has a blocking vertex $u$ or (b) $H\del\st_H(v)$ has a minimal balancing set $X$ such that $|X|\leq2$ or $H[X]$ is a triangle. When (a) happens, $u\in S$ by \ref{10.1}, so (1) holds. Assume that (b) happens. Since each cycle $C$ in $H-v$ with $|C\cap X|=1$ is unbalanced by Lemma \ref{2.1}, we have $X\subseteq \st(S)$, so $H-(S\cup X\cup\{v\})$ is 2-connected. Then $X\cup E(H[S\cup \{v\}])$ is a balancing set of $H$ with rank at most 5 by Lemma \ref{5.2} (2), a contradiction to Lemmas \ref{minimal balancing} and \ref{rank of balance set}. 

Now, we prove that (2) holds. Assume that, besides $S$, the graph $H$ has two other blocking pairs $S_1, S_2$. By \ref{10.1}, we may assume that $S_1=\{u,v\}$ and $S_2=\{u',v'\}$, where $u,u'\in S$ and $v,v'\in V(H)-S$. When $u=u'$, let $w$ be the unique vertex in $S-\{u\}$. Since $v,v', w$ are distinct blocking vertices of $H-u$, at least one pair of vertices in $\{v,v',w\}$ is a cut of $H-u$ by Lemma \ref{blocking vertex} (2), so $H$ is not 4-connected, a contradiction. Hence, $u\neq u'$, implying $v=v'$ using \ref{10.1} again.
Let $E_1$ be the set of edges from $u$ to $u'$. Since $u, u'$  are blocking vertices of $G-v$, by Lemma \ref{blocking vertex} (1) or Lemma \ref{blocking vertex} (2), either $\{u,u'\}$ is a cut of $G-v$ or $\{v\}\cup E_1$ is a balancing set. Since $H$ is 5-connected, $\{v\}\cup E_1$ is a balancing set. Since $u,v$  are blocking vertices of $G-u'$, by symmetry we have that $\{u'\}\cup E_2$ is a balancing set, where $E_2$ is the set of edges between $u$ and $v$. Hence, by Lemma \ref{5.2} (2), the set of all edges between $u, u'$ and $v$ is a balancing set, a contradiction to Lemmas \ref{minimal balancing} and \ref{rank of balance set}.
\end{proof}

\begin{lemma}\label{case<3}
Let $H$  be a $6$-connected unbalanced framework for a $7$-connected matroid $N$ with $|V(H)|\geq7$. Assume that $H$ has a blocking pair $S$ and $H-\loops(H)$ has no blocking vertices. Then $H-S$ has a balanced cycle $C$ such that all vertices in $V(C)$ are fixed in $H\del f$ and $H$ for every edge $f$ of $C$.
\end{lemma}
\begin{proof}
By Lemma \ref{case<3-1} (2), $H-S$ has a balanced cycle $C$ such that each vertex in $V(C)$ is not contained in a blocking pair of $H$. Lemma \ref{case<3-1} (1) implies that all vertices in $V(C)$ are fixed in $H$. 
Let $f$ be an arbitrary edge in $C$. Assume that the lemma does not hold for $f$. Then there is some vertex $v\in V(C)$ that is not fixed in $H\del f$. By Lemma \ref{case<3-1} (1), $\{u,v\}$ is a blocking pair of $H\del f$ for some $u\in S$. Since $\{u,v\}$ is not a blocking pair of $H$, there is a minimal balancing set $X$ of $H$ with $f\in X\subseteq\{u,v,f\}$. By Lemma \ref{5.2} (1), $X\cup S$ contains a cut of $H$, a contradiction to the fact that $H$ is 6-connected. 
\end{proof}

To prove Theorem \ref{main-result}, we need one more  result. Tutte \cite{Tutte} proved

\begin{theorem}(\cite{Oxley}, Theorem 10.3.1.)\label{ex-mi-of-bi}
A matroid is graphic if and only if it has no minor isomorphic to $U_{2,4}, F_7, F^*_7, M^*(K_5)$ and $M^*(K_{3,3})$.
\end{theorem}

Now, we prove Theorem \ref{main-result}, which is restated here in a slightly different way.

\begin{theorem}\label{main-result+1}
Let $M$ be an excluded minor for the class of quasi-graphic matroids. Then $M$ is isomorphic to $U_{3,7}$ or $U_{4,7}$, or $M$ is not $9$-connected.
\end{theorem}
\begin{proof}
Assume that $M$ is 9-connected. When $r(M)\leq 8$, it follows from Theorem \ref{rank<8} that $M$ is isomorphic to $U_{3,7}$ or $U_{4,7}$. So we may assume that $r(M)\geq9$. Since $M$ is non-graphic and the matroids in Theorem \ref{ex-mi-of-bi} each have a cocircuit of size less than 9, there is an element $e$ of $M$ such that $M\del e$ is non-graphic by Theorem \ref{ex-mi-of-bi}. Let $G$ be a 7-connected framework for $M\del e$ with $|V(G)|\geq9$.

First, consider the case that $G-\loops(G)$ has no blocking vertices and $G$ has no blocking pairs. By Lemma \ref{case>2+}, there is an edge $f$ of $G$ such that $G\del f$ has no blocking pairs and $G-(\loops(G)\cup\{f\})$ has no blocking vertices. It follows from Lemmas \ref{case>2} and \ref{almost fix} that $G\del f$ is a unique framework for $M\del e,f$. Then  $G\del f$ can be extended to a framework for $M\del f$. Moreover, since $G$ has no blocking pairs, $G\del f$ has no blocking vertices, so $M$ is quasi-graphic by Lemma \ref{2-way-extendable}, a contradiction.

Secondly, consider the case that $G$ has a blocking pair or $G-\loops(G)$ has a blocking vertex. When  $G-\loops(G)$ has a blocking vertex $v$ and $G$ is not a lifted-graphic representation for $N$, let $G'$ be obtained from $G$ by changing each loop in $\loops(G)-\loops_G(v)$ to a link joining its original end and $v$; otherwise, set $G'=G$. By Theorem \ref{ub loop}, $G'$ is also a 7-connected framework for $M\del e$ that has a blocking pair or a blocking vertex. By Lemma \ref{case:almost-balance} or Lemma \ref{case<3}, there is a balanced cycle $C$ of $G'$ such that all vertices in $V(C)$ are fixed in $G'\del f$ for each edge $f$ in $C$. Since $M\del f$ is quasi-graphic, so is $M$ by Lemma \ref{extendable}, a contradiction.
\end{proof}

\section*{Acknowledgements}
The author thanks the referees for their careful reading of this manuscript and thanks one of them for pointing out an error in the original proof of Lemma \ref{graphic matroid+1}.




\end{document}